\documentclass{article}\begin{document}
\centerline{\bf On the Evaluation of the Fifth Degree Elliptic Singular Moduli}  
\vskip .4in
\centerline{N.D.Bagis}
\centerline{Stenimahou 5 Edessa, Pellas 58200}
\centerline{Edessa, Greece}
\centerline{nikosbagis@hotmail.gr}
\vskip .2in
\[
\]
\textbf{Keywords}: Singular Moduli; Algebraic Numbers; Ramanujan; Continued Fractions; Elliptic Functions; Modular equations; Iterations; Polynomials; Pi;
\[
\]
\centerline{\bf Abstract}

\begin{quote}
We find in a algebraic radicals way the value of singular moduli $k_{25^nr_0}$ for any integer $n$ knowing only two consecutive values $k_{r_0}$ and $k_{r_0/25}$.

\end{quote}

\section{Introduction and Definitions}

Modular equations of $k_r$ (the elliptic singular moduli), have considered and have been discussed in the last 200 years by many great Mathematicians. They play very important role in several problems. The construction of $\pi$ approximation formulas, the evaluation of the famous Rogers-Ramanujan and  similar continued fractions, the solution of the quintic and sextic equation, the evaluation of the elliptic integrals in modular bases other than the classical (i.e. the cubic, the quartic and the fifth), the evaluation of the derivatives of Jacobi theta functions and many other problems of mathematics (see [1],[2],[4],[5],[6],[11],[15],[16]).\\
The only known solvable modular equations of $k_r$ where that of 2-nd and 3-rd degree. The partial solution of 5-th degree modular equation presented here is a new and important result.\\   
As application of this result we give an evaluation, in a closed form, of a quintic iteration formula for $1/\pi$, constructed by the Borwein's brothers and Bailey (see [11],[5] pg.175 and [1] pg.269).\\
We begin with the definition of the complete elliptic integral of the first kind which is (see [3],[4],[5]):
\begin{equation}
K(x)=\frac{\pi}{2}{}_2F_1\left(\frac{1}{2},\frac{1}{2};1,x^2\right)=\int^{\pi/2}_{0}\frac{1}{\sqrt{1-x^2\sin(t)^2}}dt
\end{equation}
It is known that the inverse elliptic nome (singular modulus or elliptic singular moduli), $k=k_r$, $k'^2_r=1-k^2_r$ is the real  solution of the equation:
\begin{equation}
\frac{K\left(k'_r\right)}{K(k_r)}=\sqrt{r}
\end{equation}
with $0<k_r<1$.\\
In what it follows we assume that $r \in \bf R^{*}_+ \rm$. If $r$ is positive rational then $k_r$ is algebraic. The function $k_r$ can evaluated in certain cases exactly (see [2],[5],[9]).\\ 
Continuing we define for $|q|<1$ the Ramanujan's eta function
\begin{equation}
f(-q):=\prod^{\infty}_{n=1}(1-q^n)
\end{equation}
For $\left|q\right|<1$, the Rogers Ramanujan continued fraction (RRCF) is defined as
\begin{equation}
R(q):=\frac{q^{1/5}}{1+}\frac{q^1}{1+}\frac{q^2}{1+}\frac{q^3}{1+}\cdots  
\end{equation} 
and the following relation of Ramanujan holds (see [1],[2],[8]):
\begin{equation}
\frac{1}{R^5(q)}-11-R^5(q)=\frac{f^6(-q)}{q f^6(-q^5)} 
\end{equation} 
We can write the eta function $f$ using elliptic functions. It holds
\begin{equation}
f(-q)^8=\frac{2^{8/3}}{\pi^4}q^{-1/3}(k_r)^{2/3}(k'_r)^{8/3}K(k_r)^4.
\end{equation}
Also holds (see [3] pg.488): 
\begin{equation}
f(-q^2)^6=\frac{2k_r k'_r K(k_r)^3}{\pi^3 q^{1/2}}
\end{equation}
\textbf{Theorem 1.1} (see [6],[7])
\begin{equation}
R^{-5}(q^2)-11-R^5(q^2)=\left(\frac{k_rk'_r}{ww'}\right)^2\left(\frac{w}{k_r}+\frac{w'}{k'_r}-\frac{w w'}{k_rk'_r}\right)^3 , 
\end{equation}
with $w^2=k_rk_{25r}$, $(w')^2=k'_rk'_{25r}$.  
\begin{equation}
k^6_r+k^3_r(-16+10k^2_r)w+15k^4_rw^2-20k^3_rw^3+15k^2_rw^4+k_r(10-16k^2_r)w^5+w^6=0
\end{equation}
\\
Once we know $k_r$ its relation with $w$ and hence with $k_{25r}$ is given from equation (9). Hence the problem of finding $k_{25r}$ reduces to solve the 6-th degree equation (9), which under the change of variables $w=\sqrt{k_rk_{25r}}$, $u^8=k^2_r$, $v^8=k^2_{25r}$ reduces to the 'depressed equation' (see [4] chapter 10):
\begin{equation}
u^6-v^6+5u^2v^2(u^2-v^2)+4uv(1-u^4v^4)=0
\end{equation} 
The depressed equation is also related with the problem of solution of the general quintic equation
\begin{equation}
ax^5+bx^4+cx^3+dx^2+ex+f=0,
\end{equation}
which can reduced with a Tchirnhausen transform into the Bring's form
\begin{equation}
x^5+ax+b=0.
\end{equation}   
The solution of the depressed equation is a relation of the form
\begin{equation}
k_{25r}=\Phi(k_r) .
\end{equation}
But such construction of the root of the depressed equation can not found in radicals (see [11]). Speaking clearly the equations (9) and (10) are not solvable in radicals. Hence we seek a way to reduce them. A way can found using the extra value of $k_{r/25}$. \\ 
In this paper we give a solution of the form
\begin{equation}
k_{25r}=\Phi(k_r,k_{r/25}) ,
\end{equation}
which can written more general
\begin{equation}
k_{25^nr_0}=\Phi_n(k_{r_0},k_{r_0/25}) \textrm{ , } n\in\bf Z\rm
\end{equation}
and $\Phi_n(x)$ are known algebraic constant functions which we evaluate them.

\section{State of the Main Theorem}

Our Main Theorem is\\
\\
\textbf{Main Theorem}\\
For $n=1,2,\ldots$ we have
\begin{equation}
k_{25^nr_0}=\sqrt{1/2-1/2\sqrt{1-4\left(k_{r_0}k'_{r_0}\right)^2\prod^{n}_{j=0}P^{(j)}\left[\sqrt[12]{\frac{k_{r_0}k'_{r_0}}{k_{r_0/25}k'_{r_0/25}}}\right]^{24}}}
\end{equation}
where the function $P$ is in radicals known function and is given from
\begin{equation}
P(x)=P[x]=U\left[Q\left[U^{*}\left[x\right]^6\right]^{1/6}\right]
\end{equation}
\begin{equation}
P^{(n)}(x)=(P\underbrace{\circ\ldots\circ}_{n} P)(x)\textrm{ , } P^{(0)}(x)=x . 
\end{equation}
The function $Q$ is that of (30) and $U$, $U^{*}$ are given from (33) and (34) below.

\section{The Reduction of the Evaluation of the 5th Degree Modular Equation}

We give below some Lemmas that will help us in the construction of proof and evaluation of function $P$ of Main Theorem.\\
\\
\textbf{Lemma 3.1} (see also [6])\\
Let $q=e^{-\pi\sqrt{r}}$ and $r$ real positive. We define
\begin{equation}
A_r:=\left(\frac{k'_r}{k'_{25r}}\right)^2\sqrt{\frac{k_r}{k_{25r}}}M_5(r)^{-3}
\end{equation}
Then
\begin{equation}
R(q)=\left(-\frac{11}{2}-\frac{A_r}{2}+\frac{1}{2}\sqrt{125+22A_r+A^2_r}\right)^{1/5}  
\end{equation}
where $M_5(r)$ is root of: $(5X-1)^5(1-X)=256 (k_r k'_r)^2 X$.\\
\textbf{Proof.}\\
Suppose that $N=n^2\mu$, where $n$ is positive integer and $\mu$ is positive real then holds that
\begin{equation}
K[n^2\mu]=M_n(\mu)K[\mu] 
\end{equation}
where $K[\mu]:=K(k_{\mu})$\\ 
The following equation for $M_5(r)$ is known
\begin{equation}
(5M_5(r)-1)^5(1-M_5(r))=256(k_rk'_r)^2M_5(r)
\end{equation}
Thus if we use (5),(6),(20), we get:
\begin{equation}
R^{-5}(q)-11-R^{5}(q)=\frac{f^6(-q)}{q f^6(-q^5)}=A_r=\left(\frac{k'_r}{k'_{25r}}\right)^2\sqrt{\frac{k_r}{k_{25r}}}M_5(r)^{-3}
\end{equation}
Solving with respect to $R(q)$ we get the result.\\
\\

Let now $q=e^{-\pi\sqrt{r}}$, $r>0$ and $v_r=R(q)$, it have been proved by Ramanujan that (see [1],[2],[8],[10]):
\begin{equation}
v^5_{r/25}=v_r\frac{1-2v_r+4v_r^2-3v_r^3+v_r^4}{1+3v_r+4v_r^2+2v_r^3+v_r^4}  
\end{equation} 
also from (23) is
\begin{equation}
A_r=R(q)^{-5}-11-R(q)^5=\frac{f(-q)^6}{qf(-q^5)^6} .
\end{equation}
Then from Lemma 3.1 
\begin{equation}
v_r=R(q)=S(A_r)=\left(-\frac{11}{2}-\frac{A_r}{2}+\frac{1}{2}\sqrt{125+22A_r+A^2_r}\right)^{1/5} . 
\end{equation}
Note that a $S$ function were defined from the 3-rd equality of (26).\\
From the above we get the following modular equation for $A_r$
$$
A_{r/25}=v_{r/25}^{-5}-11-v_{r/25}^{5}=F(v_r)=
$$
\begin{equation}
=\left(v_r\frac{1-2v_r+4v_r^2-3v_r^3+v_r^4}{1+3v_r+4v_r^2+2v_r^3+v_r^4}\right)^{-1}-11-\left(v_r\frac{1-2v_r+4v_r^2-3v_r^3+v_r^4}{1+3v_r+4v_r^2+2v_r^3+v_r^4}\right)
\end{equation}
and from (26) replacing $v_r$ in terms of $A_r$
\begin{equation}
A_{r/25}=F(v_r)=F\left(S\left(A_r\right)\right)=\left(F\circ S\right)(A_r)=Q(A_r) ,
\end{equation}
which after algebraic simplification with Mathematica program we get\\
\\
\textbf{Lemma 3.2}\\
If $q=e^{-\pi\sqrt{r}}$, $r>0$ and
\begin{equation}
A_r=\frac{f(-q)^6}{qf(-q^5)^6}\textrm{ , then }A_{r/25}=Q(A_r)
\end{equation}
where
\begin{equation}
Q(x)=\frac{\left(-1-e^{\frac{1}{5}y}+e^{\frac{2}{5} y}\right)^5}{\left(e^{\frac{1}{5}y}-e^{\frac{2}{5}y}+2 e^{\frac{3}{5} y}-3 e^{\frac{4}{5}y}+5 e^{y}+3 e^{\frac{6}{5}y}+2 e^{\frac{7}{5}y}+e^{\frac{8}{5} y}+e^{\frac{9}{5}y}\right)}
\end{equation}
and
\begin{equation}
y=\textrm{arcsinh}\left(\frac{11+x}{2}\right).
\end{equation}
\\

Note that inserting (31) to (30) and simplifying, we get an  algebraic function, but for simplicity and more concentrated results we leave it as it is.\\   
\\

Consider now the following equation which appear in the Lemma 3.3 below
\begin{equation}
\frac{X^2}{\sqrt{5}Y}-\frac{\sqrt{5}Y}{X^2}=\frac{1}{\sqrt{5}}\left(Y^3-Y^{-3}\right)
\end{equation}
This equation is solvable in radicals with respect to $Y$ and $X$ also. One can find the solution
\begin{equation}
Y=U(X)=\sqrt{-\frac{5}{3 X^2}+\frac{25}{3 X^2 h(X)}+\frac{X^4}{h(X)}+\frac{h(X)}{3 X^2}}
\end{equation}
where
$$
h(x)=\left(-125-9x^6+3\sqrt{3}\sqrt{-125x^6-22x^{12}-x^{18}}\right)^{1/3}
$$
The solution of (32) with respect to $X$ is
\begin{equation}
X=U^{*}(Y)=\sqrt{-\frac{1}{2 Y^2}+\frac{Y^4}{2}+\frac{\sqrt{1+18 Y^6+Y^{12}}}{2 Y^2}} .
\end{equation}
\\
\textbf{Lemma 3.3} (see [8])\\
If $G_r$ denotes the Weber class invariant 
$$
A'=\frac{f(-q^2)}{q^{1/3}f(-q^{10})}=\left(A_{4r}\right)^{1/6} \textrm{ and } V^{'}=\frac{G_{25r}}{G_r},
$$
then
\begin{equation}
\frac{A^{'2}}{\sqrt{5}V^{'}}-\frac{\sqrt{5}V^{'}}{A^{'2}}=\frac{1}{\sqrt{5}}\left(V^{'3}-V^{'-3}\right)
\end{equation}
\\
\textbf{Note.} For the Weber class invariant one can see [14],[2].\\
\\ 
We state now our first theorem\\
\\
\textbf{Theorem 3.1}\\
For the Weber class invariant $G_r$ holds
\begin{equation}
\frac{G_r}{G_{r/25}}=U\left[Q\left[U^{*}\left[\frac{G_{25r}}{G_r}\right]^6\right]^{1/6}\right]
\end{equation}
\textbf{Proof.}\\
Set $$A=\left(A_{4r/25}\right)^{1/6} \textrm{ and }V'=\frac{G_{r}}{G_{r/25}} , $$
then from Lemmas 3.2 and 3.3 and from relations (29),(30),(32),(33) and (34) we have 
$$
\frac{G_r}{G_{r/25}}=U\left[\left(A_{4r/25}\right)^{1/6}\right]=U\left[Q\left(A_{4r}\right)^{1/6}\right]=U\left[Q\left[U^{*}\left(\frac{G_{25r}}{G_{r}}\right)^6\right]^{1/6}\right]
$$
which completes the proof.  
\[
\]

Continuing we have
\begin{equation}
G_r=2^{-1/12}(k_rk'_r)^{-1/12}
\end{equation}
hence
\begin{equation}
\left(\frac{k_{r}k'_{r}}{k_{r/25}k'_{r/25}}\right)^{-1/12}=U\left[Q\left(U^{*}\left(\left(\frac{k_{25r}k'_{25r}}{k_{r}k'_{r}}\right)^{-1/12}\right)^6\right)^{1/6}\right]
\end{equation}
From the identity $$k_{1/r}=k'_r$$ we get
\begin{equation}
\left(\frac{k_{1/r}k'_{1/r}}{k_{25/r}k'_{25/r}}\right)^{-1}=U\left[Q\left(U^{*}\left(\left(\frac{k_{1/(25r)}k'_{1/(25r)}}{k_{1/r}k'_{1/r}}\right)^{-1/12}\right)^6\right)^{1/6}\right]^{12}
\end{equation}
and setting $r\rightarrow 1/r$ we lead to\\
\\
\textbf{Theorem 3.2}\\
If $r\in \bf R^{*}_{+}\rm$, then
\begin{equation}
\sqrt[12]{\frac{k_{25r}k'_{25r}}{k_{r}k'_{r}}}=U\left[Q\left[U^{*}\left[\sqrt[12]{\frac{k_{r}k'_{r}}{k_{r/25}k'_{r/25}}}\right]^6\right]^{1/6}\right]
\end{equation}
\\

Hence knowing $k_{r_0}$ and $k_{r_0/25}$ we can evaluate in closed form the $k_{25r_0}$. If we repeat the process we can find any higher or lower order of $k_{25^nr_0}$ in closed radicals form, when $n\in\bf Z\rm$ and easily get (16).\\ 
Observe that a similar formula to (16) for the evaluation of $k_{r_0/25^n}$, $n=1,2,\ldots$ can extracted from (39).\\
\\ 
\textbf{Example 1.}
$$
k_{1/5}=\sqrt{\frac{9+4 \sqrt{5}+2 \sqrt{38+17 \sqrt{5}}}{18+8 \sqrt{5}}}
$$
$$
k_{5}=\sqrt{\frac{9+4 \sqrt{5}-2 \sqrt{38+17 \sqrt{5}}}{18+8 \sqrt{5}}}
$$
\begin{equation}
k_{125}=\sqrt{\frac{1}{2}-\frac{1}{2}\sqrt{1-(9-4\sqrt{5})P[1]^2}}  
\end{equation}
where $P(x)=P[x]$ is that of (17) and $P[1]$ is a algebraic radical number which has too complicated form to present it here. Its minimal polynomial is 
$$
1-5 x^2-10x^3-15 x^4-22 x^5-15 x^6-10x^7-5 x^8+x^{10}=0
$$
\\
The same holds also for the next\\
\\
\textbf{Example 2.}\\
for $r_0=25$ it is
$$
k_{r_0/25}=k_1=\frac{1}{\sqrt{2}}
$$
and
$$
k_{r_0}=k_{25}=\frac{1}{\sqrt{2 \left(51841+23184 \sqrt{5}+12 \sqrt{37325880+16692641 \sqrt{5}}\right)}}
$$
Hence
\begin{equation}
k_{625}=\sqrt{\frac{1}{2}-\frac{1}{2}\sqrt{1-\left(51841-23184\sqrt{5}\right)P\left[\frac{\sqrt{5}-1	}{2}\right]^{24}}}
\end{equation}
\[
\]
Hence from 2-nd and 3-rd degree modular equations (see [1] chapter 19), we can evaluate every $k_{r}$ which $r$ is of the form $r=4^k9^l25^nr_0$, when $k_{r_0}$ and $k_{r_0/25}$ are known and $k,l,n\in\bf Z\rm$.

\section{The fifth degree singular moduli and approximations to $1/\pi$}

In [5] pg.175 J.M. Borwein and P.B. Borwein consider the following sharp convergent approximation algorithm to $\pi$ (see also [11]):\\
Consider $\alpha_0:=\alpha(r_0)$ where $\alpha(n)$ is the elliptic alpha function and $u_0=\sqrt[8]{1-k_{r_0}^2}$ set also $u_{n}=(k_{25^{n-1}r_0})^{1/4}$, which now are given from (16) and are in closed form radicals. Using this way we are able to constuct approximations not depending on numerical values of singular moduli $k_r$, but from exact values.\\
Let 
$$
x_n:=2u_nu^5_{n+1}\textrm{, }y_n:=2u_{n}^5u_{n+1}
$$ 
$$
a_n:=u_{n}^2+5u^2_{n+1}+2x_n\textrm{, }b_n:=5u_{n}^2+u^2_{n+1}-2y_n\textrm{, }\gamma_n:=\frac{a_n}{b_n}
$$
finally
$$
\delta_n:=4^{-1}a^{-1}_n(1-u^8_{n+1})\left[5(u^2_{n+1}+x_n)+\gamma_{n}(y_n-u^2_{n+1})\right]+
$$
$$
+4^{-1}b_n^{-1}(1-u^8_n)[u^2_{n}+x_n+5\gamma_{n}(y_n-u^2_{n})]
$$
Then
\begin{equation}
\alpha_{n+1}:=5\gamma_n\alpha_n+5^{n+1}\sqrt{r_0}(\delta_n+u^8_{n+1}-\gamma_nu^8_n)
\end{equation}
and 
\begin{equation}
0<\alpha_n-\pi^{-1}<16\cdot5^n\sqrt{r_0}e^{-5^n\sqrt{r_0}\pi}
\end{equation}
for $r_05^{2n}\geq 1$.\\
Hence for every initial condition given to $r_0>0$, with $k_{r_0}$ and $k_{r_0/25}$ known, we lead to a closed form iteration formula approximating $1/\pi$. Actually the thirteen iterations of the above algorithm give the first one billion digits of $\pi$. Note also that algorithms presented in [11],[13] are for a specific initial value and do not cover all values of $r_0$. Here if we know only $k_{r_0}$ and $k_{r_0/25}$ we can find all $u_n$ and construct the iteration for $1/\pi$.\\
Our idea can generalized also for the 10th degree modular equation of $k_r$. The 10th degree modular equation of Rogers Ramanujan continued fraction is solvable and also can put in the form $v_{r/100}=\phi(v_r)$ (see [12]). By this way one can solve with initial conditions $k_{r_0}$ and $k_{r_0/100}$ the general 10-th degree modular equation of $k_r$, if finds an analogue of Lemma 3.3. But it is not so economical since the modular equation of 10-th degree can be splitten to that of 2-nd and 5-th degree.

\[
\]

\centerline{\bf References}\vskip .2in

\noindent

[1]: B.C. Berndt. 'Ramanujan`s Notebooks Part III'. Springer Verlag, New York (1991)

[2]: B.C. Berndt. 'Ramanujan`s Notebooks Part V'. Springer Verlag, New York (1998)

[3]: E.T. Whittaker and G.N. Watson. 'A course on Modern Analysis'. Cambridge U.P. (1927)

[4]: J.V. Armitage W.F. Eberlein. 'Elliptic Functions'. 
Cambridge University Press. (2006)

[5]: J.M. Borwein and P.B. Borwein. 'Pi and the AGM'. John Wiley and Sons, Inc. New York, Chichester, Brisbane, Toronto, Singapore. (1987) 

[6]: N.D. Bagis. 'The complete evaluation of Rogers-Ramanujan and other continued fractions with elliptic functions'. arXiv:1008.1304v1 [math.GM] 7 Aug 2010

[7]: N.D. Bagis. 'Parametric Evaluations of the Rogers Ramanujan Continued Fraction'. International Journal of Mathematics and Mathematical Sciences. Vol. 2011

[8]: B.C. Berndt, H.H. Chan, S.S. Huang, S.Y. Kang, J. Sohn and S.H. Son. 'The Rogers-Ramanujan Continued Fraction'. J. Comput. Appl. Math. 105 (1999), 9-24

[9]: D. Broadhurst. 'Solutions by radicals at Singular Values $k_N$ from New Class Invariants for $N\equiv3\;\; mod\;\; 8$'. arXiv:0807.2976 (math-phy)

[10]: W. Duke. 'Continued Fractions and Modular Functions'. Bull. Amer. Math. Soc.(N.S.) 42 (2005), 137-162.
    
[11]: J.M. Borwein, P.B. Borwein and D.H. Bailey. 'Ramanujan, Modular Equations, and Approximations to Pi or How to Compute One Billion Digits of Pi'. Amer. Math. Monthly (96)(1989), 201-219

[12]: M. Trott. 'Modular Equations of the Rogers-Ramanujan Continued Fraction'. Mathematica Journal. 9,(2004), 314-333

[13]: H.H. Chan, S. Cooper, and W.C. Liaw. 'The Rogers-Ramanujan Continued Fraction and a Quintic Iteration for 1/$\pi$'. Proc. of the Amer. Math. Soc. Vol.135, No.11, (2007), 3417-3424 

[14]: B.C. Berndt and H.H. Chan. 'Notes on Ramanujan's Singular Moduli'. The Proceedings of the Fifth Conference of the Canadian Number Theory Association, edited by R.Gupta and K.S. Williams, (1998), 7-16

[15]: N.D. Bagis and M.L. Glasser. 'Conjectures on the evaluation of alternative modular bases and formulas approximating $1/\pi$'. Journal of Number Theory. (Elsevier), (2012).

[16]: N.D. Bagis. 'A General Method for Constructing Ramanujan-Type Formulas for Powers of $1/\pi$'. The Mathematica Journal, (2013)

\end{document}